\documentclass[reqno,12pt]{amsart}
\usepackage{a4wide}
\usepackage{amsmath}
\usepackage{amssymb}
\usepackage{amsthm}
\usepackage{mathrsfs}
\usepackage{graphicx}

\numberwithin{equation}{section}

\def\accentsfrancais{applemac}
   \RequirePackage[\accentsfrancais]{inputenc}

\newtheorem{theorem}{Theorem}

\newtheorem{proposition}{Proposition}

\theoremstyle{remark}

\def\R{\mathbb{R}}

\def\gam{\gamma}
\def\Pref{{\rm Pref}}

\def\N{\mathbb{N}}
\def\P{\mathbb{P}}
\def\R{\mathbb{R}}

\def\Z{\mathbb{Z}}

\def\Om{\Omega}
\def\({\left(}
\def\){\right)}
\def\[{\left[}
\def\]{\right]}

\def\mk{\medskip}
\def\sk{\smallskip}

\begin{document}

\title[Some Aspects of Multifractal analysis]
{Some Aspects of Multifractal analysis}

\author[]{Ai-Hua Fan}

\address{Ai Hua FAN, LAMFA, UMR 7352,  CNRS, Universit\'e de Picardie, 33 rue
Saint Leu, 80039 Amiens, France.
E-mail: ai-hua.fan@u-picardie.fr}



\begin{abstract}
The aim of this survey  is to present some aspects of multifractal analysis around the recently developed subject   of multiple ergodic averages. Related topics include  dimensions of measures, oriented walks, Riesz products etc. The exposition on the multifractal analysis of  multiple ergodic averages is mainly based on \cite{FLM,KPS,FSW}.
\medskip

{\em This paper is published in Geometry and Analysis of Fractals (pp 115-145), Ed. D. J. Feng and K. S. Lau,  Springer Proceedings in Mathematics and Statistics Vol. 88, Springer (2014).}

\end{abstract}

\maketitle

\section{Introduction}
\label{sec:introduction}

Multifractal problems can be put into the following frame.
Let $(X, d)$ be a metric space and $\mathcal{P}(x)$ a property (quantitative or qualitative) depending on a point $x$ of the space $X$.
For any prescribed  property $\mathbf{P}$, we look at the set of those points $x$ which have the property $\mathbf{P}$:
$$
   E(\mathbf{P}) =\{x \in X: \mathcal{P}(x) = \mathbf{P}\}.
$$
The size of the sets $E(\mathbf{P})$ for different $\mathbf{P}$'s is problematic in the multifractal analysis. According to popular folklore, the function $x \mapsto \mathcal{P}(x)$ is multifractal
if $E(\mathbf{P})$ is not empty for  uncountably many properties $\mathbf{P}$. Usually the size
of a set $A$ in $X$ is described by its Hausdorff dimension $\dim_H A$ or its packing dimension
$\dim_P A$, or its topological entropy in a dynamical setting. See \cite{Falc,Mattila1995} for the dimension theory and \cite{Bowen73,Pesin} for the notion of topological entropy.

Seeds of multifractals were sown in B. Mandelbrot's works on multiplicative chaos
in 1970's \cite{Mandelbrot1974a,Mandelbrot1974b}. 
First rigorous results are due to J. P. Kahane and J. Peyr\`ere \cite{KP1976}. The concept of multifractality came from geophysics and theoretical physics.  At the beginning,
U. Frisch and G. Parisi~\cite{FrishParisi},  H. G. Hentschel and I. Procaccia~\cite{hentchel}
had the rather vague idea of mixture of subsets of different dimensions each of which has
a given H\"{o}lder singularity exponent. The multifractal formalism  became clearer in the 1980-90's
in the works of  T. C. Halsey-M. H. Jensen-L.P. Kadanoff-I. Procsccia-B. I. Shraiman \cite{Halsey}, of P. Collet-J. L. Lebowitz-A. Porzio~\cite{collet}, and of G. Grown-G. Michon-J. Peyr\`ere \cite{BMP}.
The multifractal formalism is tightly related to the thermodynamics and D. Ruelle \cite{Ruelle} was the
first to use the thermodynamical formalism to compute the Hausdorff dimensions of some Julia sets.

The research on the subject has been very active and very fruitful since four decades.
The first most studied multifractal quantity is the local dimension of a Borel measure $\mu$ on $X$. Recall that
the local lower dimension of $\mu$ at $x\in X$ is defined by
$$
    \underline{D}_\mu(x) = \liminf_{r\to 0} \frac{\log \mu(B(x, r))}{\log r}
$$
where $B(x, r)$ denotes the ball centered at $x$ with radius $r$. The local upper
dimension  $\overline{D}_\mu(x)$ is similarly defined.
There is a huge literature on this subject. 
Let us mention another example, H\"{o}lder exponent of a function (see \cite{Jafsiam}).
Let $\alpha >0$. Let $f: \mathbb{R}^d \to \mathbb{R}$ be a function and  $x\in
\mathbb{R}^d$ be a fixed point. We say $f$ is $\alpha$-H\"{o}lder at
$x$ and we write $f\in C^\alpha(x)$ if there exist two constants $
\delta >0$ and $C>0$ and a polynomial $P(x)$ of degree strictly smaller than $\alpha$  such that
$$
   |y-x|<\delta \Rightarrow |f(y) -f(x) - P(y-x)|\le C |y-x|^\alpha.
$$
The H\"{o}lder exponent of $f$ at $x$ is then defined by
$$
   h_f(x) =\sup\{\alpha>0: f\in C^\alpha(x)\}.
$$

The multifractal analysis has now become  a set of tools applicable in analysis, probability
and stochastic processus, number
theory, ergodic theory and dynamical systems etc if we don't account applications in physics
and other sciences.
\mk

The main goal of  this paper is to present some problems in the multifractal analysis of
Birkhoff ergodic averages, especially of  multiple Birkhoff ergodic averages, and some related topics
like dimensions of measures and Riesz products as tools, and oriented walks as similar subject.

\mk

Let $T: X \to X$ be a map from $X$ into $X$. We consider the dynamical system
$(X,T)$. The main concern about the system
is the behavior of an orbit $\{T^n x\}$ of a given point $x\in X$. Some aspects of the behavior of the orbit
may be described by  the so-called Birkhoff averages
\begin{equation}
\label{def_Sn}
 A_nf (x) =  \frac{1}{n} \sum_{k=0}^{n-1} f(T^k x)
\end{equation}
where $f: \mathbb{R}^d \to \mathbb{R}$ is a given function, called observable. We refer to \cite{Walters1982}
for basic facts in the theory of dynamical system and ergodic  theory.

\mk

The famous and fundamental
Birkhoff Ergodic Theorem states that for any $T$-invariant ergodic probability measure
$\mu$ on $X$ and for any integrable function $ f\in L^1(\mu)$, the limit
$$ \lim_{n\to +\infty} \ A_nf(x) = \int_X f \, d\mu
$$
holds for $\mu$-almost all points $x\in X$.  Even if $\mu$ is only $T$-invariant but not ergodic, the limit still exists for $\mu$-almost all points $x\in X$. For many dynamical systems, there is a rich class
of invariant measures and so that the limit of Birkhoff averages $A_nf(x)$ may vary for different points $x$. This variety reflects the chaotic feature of the dynamical system. A typical example is the
doubling dynamics $x \mapsto 2 x (\!\!\!\mod 1)$ on the unit circle $\mathbb{T}:=\mathbb{R}/\mathbb{Z}$.

The multifractal analysis of Birkhoff ergodic averages provides a way to study the chaotic feature
of the dynamics. Let $\alpha\in \R$. We define the $\alpha$-level set
$$E_f(\alpha) = \left\{ x\in X: \lim_{n\to +\infty} A_nf(x)   = \alpha \right\}.$$
The purpose of the multifractal analysis is the determination of the size of  sets $ E_f(\alpha)$.
If the Hausdorff dimension is used as measuring device, we are led to the Haudorff multifractal
spectrum of $f$:
$$
   \mathbb{R} \ni \alpha \mapsto  d_H(\alpha): = \dim_H E_f(\alpha).
$$
The existing works show that in many cases it is possible to compute the spectrum $d_H(\cdot)$
and it is also possible to distinguish a nice invariant measure sitting on the set $F_f(\alpha)$
for each $\alpha$. By ``nice" we mean that the measure is supported by $E_f(\alpha)$ and its dimension is equal to that of  $E_f(\alpha)$. The dimension of a measure is defined to be the dimension of the ``smallest" Borel support of the measure. Therefore the nice measure is a maximal measure in the sense that it attains the maximum among all measures supported by $E_f(\alpha)$. This maximal measure may be
invariant, ergodic and even mixing. Some other nice properties are also shared by the maximal measure.
For this well studied  classic ergodic averages, see for example \cite{BSS,FF,FFW,FLP2008,FS2003,FengLauWu,Olivier,TV}. Let us also mention some useful tools
for dimension estimation \cite{BarralSeuret,BV,Durand,FST}.

\sk
The multiple ergodic theory started almost at the same time of the development of multifractal analysis. It started with F\"{u}rstenberg's proof of Szemer\'edi theorem on the existence of
arbitrary long arithmetic sequence in a set of integers of positive density \cite{Furstenberg}. This theory involves several dynamics rather than one dynamics. Let $T_1$, $T_2$, ..., $T_d$ be $d$ transformations on a space $X$.
We assume that they are commuting each other and preserving  a given probability measure $\mu$.
For $d$ measurable functions $F=(f_1, \cdots, f_d)$ we define the  multiple Birkhoff averages
by
 \begin{equation}
\label{def_generalergaverage-0}
 A_n F(x) :=\frac{1}{n} \sum_{k=0}^{n-1} f_1(T_1^k x)f_2( T_2^k  x)\cdots f_d(T_d^k x).
\end{equation}
An inspiring example is the couple $(\tau_2, \tau_3)$ on the circle $\mathbb{T}$ where
$$
\tau_2 x = 2x \ (\!\!\!\!\!\mod 1), \qquad\tau_3 x = 3x \ (\!\!\!\!\!\mod 1).
$$

The mixture of the dynamics $T_1, T_2, \cdots, T_d$ is much more difficult to understand.
After the first works of F\"{u}rstenberg-Weiss \cite{FurstenbergWeiss} and of Conze-Lesigne \cite{ConzeLesigne},
B. Host and B. Kra \cite{HK1} proved  the $L^2$-convergence of $A_nF$ when $T_j=T^j$
(powers of a fixed dynamics) and $f_j \in L^\infty(\mu)$. For the almost every convergence,
results are sparse. J. Bourgain proved the almost everywhere convergence when $d=2$.
We should point out that even if the limit exists, it is not easy to recognize the limit.
In particular, the limit may not be constant for some ergodic measures. For nilsystems, explicit formula for the limit was known to E. Lesigne \cite{Lesigne} and T. Ziegler \cite{Ziegler}. Anyway, not like the "simple"
ergodic theory, the multiple ergodic
theory has not yet reached its maturity.

Although the multiple ergodic
theory has been still developing, this situation doesn't prevent us from investigating the
multifractal feature of multiple systems.
Let us  consider the following general set-up.
For a given observable $\Phi:X^d \to \R$,  we consider the Multiple Ergodic Averages
\begin{equation}
\label{def_generalergaverage}
 A_n \Phi(x) :=\frac{1}{n} \sum_{k=0}^{n-1} \Phi(T_1^k x, T_2^k  x, ..., T_d^k x).
\end{equation}
The Birkhoff averages (\ref{def_generalergaverage-0}) correspond to
the special case of tensor product $\Phi=f_1\otimes f_2\otimes \cdots \otimes f_d$.

It is natural to introduce the following generalization of multifractal
spectrum. In  this paper, we will only consider the case
where $T_j = T^j$ for $1\le j \le d$. So
\begin{equation}
\label{def_ergaverage}
A_n \Phi(x)=\frac{1}{n} \sum_{k=0}^{n-1} \Phi(T^k x, T^{2k}  x, ..., T^{dk} x).
\end{equation}
The multiple Hausdorff spectrum of the observable $\Phi$ is defined by
$$
    d_H(\alpha) = \dim_H E_\Phi(\alpha), \qquad (\alpha \in \mathbb{R})
$$
where
\begin{equation}
\label{def_Ealpha}
E_\Phi(\alpha)=\left\{ x\in X:  \lim_{n\to +\infty} \ \frac{1}{n} \sum_{k=0}^{n-1} \Phi(T^k x, T^{2k}  x, ..., T^{dk} x) =\alpha\right\}.
\end{equation}

As we said, the classical theory ($d=1$) is well developed.
When $d>1$, there are  several results on the multifractal analysis of the limit of the averages $A_n\Phi$ in some special cases; but questions remain
largely unanswered.  In Section 2, we will present the first results obtained in \cite{FLM} in a very special case which
give us a feeling of the problem and illustrate the difficulty of the problem. One result is  the Hausdorff spectrum  obtained   by using Riesz products, a tool borrowed from Fourier analysis. Another  result concerns the box dimension of a multiplicatively invariant set. These two kinds of result
are respectively generalized in \cite{FSW} and \cite{KPS,PSSS} in some general setting.
In Section 5, we will present these results.
As we mentioned,
the local dimension of a measure was the first study object of multifractal analysis. In Section 3,
we will give an account of dimensions of measures which are related to the local dimension and have their own interests.  The idea of using Riesz products is inspired by a work on oriented walks \cite{Fan2000}
to which Section 4 will be devoted. 
In the last section,  we will collect  some remarks and open problems. 

\section{First multifractal results on the multiple ergodic averages }

The question of computing the dimension of  $E_\Phi(\alpha)$  was  raised by Fan, Liao and Ma in \cite{FLM}, where the following special case was studied: $X=: \mathbb{M}_2:=\{-1, 1\}^\mathbb{N^*}$,
$T$ is the shift, $d\ge 1$  and $\Phi$ is the function
$$
    \Phi(x^{(1)}, x^{(2)}, \cdots, x^{(d)})= x^{(1)}_1 x^{(2)}_1\cdots x^{(d)}_1, \quad (
    x^{(1)}, x^{(2)}, \cdots, x^{(d)})\in \mathbb{M}_2^d.
$$
Note that $x_1^{(j)}$ is the first coordinate of $x^{(j)}$.
We consider $\mathbb{M}_2$ as the infinite product of the multiplicative group $\{-1, 1\}$. Then
the function $x\mapsto x_j$ is a group character of $\mathbb{M}_2$, called Rademacher function and
$$
    \Phi( T^k x, T^{2d}x, \cdots, T^{dk}x)
    = x_k x_{2k}\cdots x_{dk}
$$
is also a group character,  called Walsh function.
Recall that $\mathbb{M}_2=\{-1,1\}^{\N^*}$, considered as a symbolic space, is  endowed with the metric
$$d(x,y) = 2 ^{-\min\{k: x_k\neq y_k\}}, \mbox{ for $x,y \in \mathbb{M}_2$.}$$

\subsection{
Multifractal spectrum of a sequence of Walsh functions} \
Under the above assumption, we have the following result.

\begin{theorem}{\rm (\cite{FLM})}\label{Thm_FLM}
For every $\alpha\in [-1,1]$, the set
$$B_\alpha := \left\{x\in \mathbb{M}_2: \lim_{n\to +\infty} \frac 1 n \sum_{k=1}^n x_kx_{2k} \cdots x_{dk} = \alpha\right\}$$
has the Hausdorff dimension
$$\dim_H B_\alpha = 1 - \frac 1 d + \frac{1} {d\log 2} H \left( \frac{1+\alpha}{2} \right),$$
where $H(t):= -t\log t -(1-t) \log(1-t)$.
\end{theorem}

A key observation is that all these
Walsh functions constitute a dissociated system in the sense of Hewitt-Zuckermann \cite{HZ}. This
allows us to define probability measures called Riesz products on the group $\mathbb{M}_2$:
$$
     \mu_b := \prod_{k=1}^\infty (1+ b x_k x_{2k}\cdots x_{dk})
     := w^*\!-\!\lim_{N\to\infty} \prod_{k=1}^N (1+ b x_k x_{2k}\cdots x_{dk} )dx
$$
for $b \in [-1, 1]$, where $dx$ denote the Haar measure on $\mathbb{M}_2$. Fortunately, the Riesz product
$\mu_\alpha$ ($b=\alpha$) is a maximizing measure on $B_\alpha$. So,  by computing
the dimension of the measure $\mu_\alpha$, we get the stated formula. We will go back to Riesz products
in Section~\ref{Sect_Dim} and to dimensions of measures in Section~\ref{Sect_OW}.

Note that the case $d=1$ is nothing but
  the well-known Besicovich-Eggleston theorem dated back to 1940's, which would be considered as the
   first result of multifractal analysis.


\mk

\subsection{Box dimension of some multiplicatively invariant set}\label{X_2}
The very motivation of \cite {FLM} was the multiple ergodic averages in the following case where
$X: = \mathbb{D}_2:=\{0,1\}^{\N^*}$, $T$ is the shift on $\mathbb{D}_2$ and $\Phi(x,y) = x_1y_1$ for $x=x_1x_2... $ and $ y=y_1y_2...\in \mathbb{D}_2$.
The space $\mathbb{D}_2$ may also be considered as the infinite product group of $\mathbb{Z}/2\mathbb{Z}$. But the function $x \mapsto x_1$ is no longer group character and the Fourier
method fails. Then the authors of \cite{FLM} proposed to look at a subset of the $0$-level set
 \begin{equation}
\label{def:E0}
E_\Phi(0) 
= \left\{x\in \Sigma_2: \lim_{n\to +\infty} \frac 1 n \sum_{k=1}^n x_k x_{2k} =0\right\}.
\end{equation}
The proposed subset is
\begin{equation}
\label{def:X2}
X_2 
= \{x\in \Sigma_2: \ \forall \, k\geq 1, \ x_k x_{2k}=0\}.
\end{equation}
This set $X_2$ has a nicer structure than $E_\Phi(0)$.
The condition $x_k x_{2k}=0$ is  imposed to all integers $k$ without exception for all points $x$ in $X_2$, while the same condition is imposed to  "most" integers $k$ for points $x$ in $E_\Phi(0)$.

\begin{theorem}\label{box-dim} {\rm (\cite{FLM})}
\label{prop_FLM}
The box dimension of $X_2$ is equal to
$$ \dim_B X_2=\frac{1}{2\log 2} \sum_{n=1}^{+\infty} \frac{ \log F_n}{2^n} ,$$
where $F_n$ is the Fibonacci sequence: $F_0=1$, $F_1=2$ and $F_{n+2}=F_{n+1}+F_n$ for $n\geq 0$.
\end{theorem}

The key idea to prove the formula  is the following observation, which is
 also one of the key points for all the obtained results up to now in different cases.
Look at the definition  \eqref{def:X2} of $X_2$. The value of the digit $x_1$ of an element $x=(x_k)\in X_2$ has an impact on the value of $x_2$, which in turn on the value of $x_4$, ... and so forth on the values of $x_{2^k}$ for all $k\geq 1$. But it has no influence on $x_3$, $x_5$, ... .
  Similarly, the value of $x_i$ for an odd integer $i$ only has influence on $x_{ i 2^k}$.  This suggests us
  the following partition
$$\N^* = \bigsqcup _{i \mbox{\rm \ \small odd}}  \   \Lambda_i, \ \  \mbox { with } \Lambda_i : =\{i\,2^n: n\geq 0\}.$$
We could say that the defining conditions of $X_2$ restricted to different $\Lambda_i$  are independent. We are then led to investigate, for each odd number $i$, the restriction of $x$ to $\Lambda_i$
which will be denoted by
$$ x_{|\Lambda_i} = x_{i}x_{i2}x_{i2^2} \ldots x_{i2^n} \ldots.$$
If we   rewrite $x_{|\Lambda_i} = z_1z_2 ...$, which is considered as a point in $\mathbb{D}_2$,
then $(z_n)$ belongs to the subshift of finite type subjected to $z_k z_{k+1} =0$.

 It is clear that
$$\dim_B X_2=\lim_{n\to\infty}\frac{\log_2N_n}{n}$$
 where $N_n$ is the cardinality of the set
$$
 \{(x_1x_2\cdots x_{n}):  x_{k}x_{2{k}}=0 \ \text{ for } k\geq 1 \text{ such that }   \ 2{k}\leq
  n\}.
$$
 Let us decompose the set of the first $n$ integers as follows
$\{1, \cdots, n\}=C_0 \sqcup C_1\sqcup
\cdots\sqcup C_m$ with
\begin{align*} C_0:&=\left\{1,  \ \ \ \ \ \ 3,\ \ \ \ \ \, 5, \ \ \ \ \dots, \  2
n_0-1 \right\}, \\
C_1:&=\left\{1\cdot 2, \ \ 3\cdot 2, \ \ 5\cdot 2, \ \dots, \
2\cdot\big(2n_1-1\big) \right\}, \\
&\dots\\
C_k:&=\left\{1\cdot 2^k , \ 3\cdot 2^k, \ 5\cdot 2^k, \ \dots, \
2^k\cdot (2n_k-1) \right\}, \\
&\dots\\
C_m:&=\left\{1\cdot 2^m \right\}.
\end{align*}
These finite sequences have different length $n_k$ ($0\le k\le m$). Actually
$n_k$ is the biggest integer such that $2^k(2n_k-1)\leq n$, i.e.
$
n_k=\left\lfloor\frac{n}{2^{k+1}}+\frac{1}{2}\right\rfloor.
$ The number $m$  is the biggest integer such that
$2^m\leq n$, i.e. $m=\lfloor\log_2 n\rfloor.$
It clear that $n_0>n_1>\cdots > n_{m-1}>n_m=1$. The
conditions $x_{\ell}x_{2\ell}=0$ with $\ell$ in different columns in
the table defining $C_0, \cdots, C_m$ are independent.
This independence allows us to count the number of possible choices for
$(x_1,\cdots, x_n)$ by the multiplication principle. First,
we have $n_m(=1)$ column  which has $m+1$ elements. We have
$F_{m+1}$ choices for those $x_{\ell}$ with $\ell$ in the first column
because $x_{\ell}x_{2\ell}$ is conditioned to be different from  the word $11$.
We repeat this argument for other columns.
 Each of the next
$n_{m-1}-n_{m}$ columns has $m$ elements, so we have $F_{m}^{n_{m-1}-n_m}$
choices for those $x_{\ell}$ with $\ell$ in these columns. By
induction, we get
$$N_n=F_{m+1}^{n_m}F_{m}^{n_{m-1}-n_m}F_{m-1}^{n_{m-2}-n_{m-1}}\cdots F_1^{n_0-n_1}.$$
To finish the computation we need to note that $\frac{n_k}{n}$ tends to $2^{-(k+1)}$ as $n$ tends to the infinity.

 The set $X_2$ is not invariant under the shift. But as observed by Kenyon, Peres and Solomyak \cite{KPS}, it is multiplicatively invariant in the sense that $M_r X_2 \subset X_2$ for all integers $r\ge 1$
 where
 $$
     M_r((x_n)) = (x_{rn}).
 $$
 The Hausdorff dimension of $X_2$ was obtained in \cite{KPS} and the gauge function
 of $X_2$ was obtained in \cite{PS2}.
 We will discuss the work done in \cite{KPS} in Section~\ref{Sect_MBA}.



\section{Dimensions of measures}\label{Sect_Dim}

Multifractal properties were first investigated for measures. The multifractal analysis of a measure is the
analysis
of the local dimension on the whole space, while the dimensions of a measure
concern with what happens on a Borel support of the measure. In \cite{Fan1989c},
lower and upper Hausdorff dimensions of a measure were introduced and systematically studied,
inspired by \cite{Peyriere1975} and \cite{Kahane1987b}. The lower and upper packing dimensions were
later studied independently in \cite{Tamashiro1995} and \cite{He1}.
Some aspects were also considered in \cite{Haase}.
A fundamental theorem in the theory of fractals is Frostman theorem. 
Howroyd \cite{Howroyd1994} and Kaufman
\cite{Kaufman} generalized it from Euclidean spaces to complete separable metric spaces.
This fundamental theorem allows us to employ the potential theory.

\subsection{Potential theory}
Let $(X, d)$ be a complete separable metric space, called Polish space.
Let $0<\alpha<\infty$.
For any locally finite Borel measure $\mu$ on $X$, we define its
  potential of order $\alpha$ by
$$
U_{\alpha}^{\mu} (x) := \int_X \frac{\text{d}\mu(y)}{(\text{d}(x,
  y))^\alpha}\quad  \,\,\,\, (x \in X)
$$
and its energy of order $\alpha$ by
$$
I_{\alpha}^{\mu} := \int_X U_{\alpha}^{\mu}(x) \text{d}\mu(x) = \int_X
\int_X \frac{\text{d}\mu(x) \text{d}\mu(y)}{(\text{d}(x, y))^\alpha}.
$$
The capacity of order $\alpha$ of a compact set $K$ in $X$
is defined by
$$
{\rm Cap}_\alpha K = \left( \inf_{\mu \in \mathcal{M}^+_1(K)}
  I_{\alpha}^{\mu} \right)^{-1}.
$$
For an arbitrary set $E$ of $X$, we define its capacity of order $\alpha$ by
$$
{\rm Cap}_\alpha E = \sup\{ {\rm Cap}_\alpha K: K \, \mbox{compact\
  contained\ in }\ E\} .
$$

For a set $E$ in $X$, we define its capacity dimension by
$$
\dim_C E = \inf \{\alpha>0: {\rm Cap}_\alpha (E) = 0\} = \sup
\{\alpha>0: {\rm Cap}_\alpha (E) > 0\}.
$$
The following is the theorem of Frostman-Kaufman-Howroyd. Frostman initially
dealt with the Euclidean space. Kaufman proved the result by generalizing a min-max
theorem on quadratic function to the mutual potential energy functional, while Howroyd
used the technique of weighted Hausdorff measures.

\begin{theorem}[\cite{Kaufman,Howroyd1994}] Let $(X, d)$ be a complete metric space.
  For any  Borel $E\subset X$, we have
$
\dim_C E = \dim_H E.
$
\end{theorem}

\subsection{Hausdorff dimensions of a measures}
An important tool to study a measure is its dimensions, which attempt to estimate the size
of the "supports" of the measure.
The idea  finds its origin in J. Peyri\`ere's works on Riesz products
\cite{Peyriere1975} and also in that of J.-P. Kahane on Dvoretzky covering \cite{Kahane1987b}. The following definitions were introduced in \cite{Fan1989c} (see also \cite{Fan1994, Falc}).

Let $(X, d)$ be a complete separable metric space. Let $\mu$ be a Borel measure on $X$.
The {\em lower Hausdorff dimension} and the {\em upper Hausdorff dimension} of a measure $\mu$  are respectively defined by
\begin{eqnarray*}
   \dim_* \mu = \inf \{ \dim_H A: \mu(A)\ >0\},\qquad \dim^* \mu  = \inf \{ \dim_H A: \mu(A^c)=0\}.
\end{eqnarray*}
It is evident  that
$
   \dim_* \mu \le \dim^* \mu.
$
When the equality holds, $\mu$ is said to be {\em unidimensional} or $\alpha$-{\em dimensional}
where $\alpha$ is the common value of $\dim_* \mu$ and $\dim^* \mu$.
The Hausdorff dimensions $\dim_* \mu$ and $\dim^* \mu$ are described by the lower local dimension
function $\underline{D}_\mu(x)$ in the following way.

\begin{theorem}[\cite{Fan1989c, Fan1994}]\label{HDM}
$$   \dim_* \mu = {\mbox{\rm ess} \inf}_\mu \ \underline{D}_\mu(x), \quad
   \dim^* \mu  = {\mbox{\rm ess} \sup}_\mu \ \underline{D}_\mu(x).
   $$
\end{theorem}
There are also a continuity-singularity criterion using Hausdorff measures
and a energy-potential criterion. Sometimes these criteria are more practical.
\begin{theorem}[\cite{Fan1989c, Fan1994}]
\begin{eqnarray*}
   \dim_* \mu &=&  \sup \{ \alpha>0: \mu \ll H^\alpha\}
               = \sup 
                       \{ \alpha>0: \mu=\sum \mu_k, I_\alpha^{\mu_k}<\infty
                       \}\\
   \dim^* \mu & =& \inf \{ \alpha>0: \ \mu \perp H^\alpha \}
               \, = \inf \, 
                        \{ \alpha>0:  U_\alpha^{\mu}(x) =\infty, \ \mu\!-\!\mbox{\rm p.p.}
                        \}.
\end{eqnarray*}
\end{theorem}

Theorem~\ref{HDM} holds for lower and upper packing dimensions of a measures,
similarly defined, if we replace  $\underline{D}(\mu, x)$ by $\overline{D}(\mu, x)$
(see \cite{Tamashiro1995} and \cite{He1}). They are denoted by $\mbox{\rm Dim}_*\mu$
and $\mbox{\rm Dim}^*\mu$.

We say $\mu$ is exact if $\dim_*\mu=\dim^*\mu =\mbox{\rm Dim}_*\mu =\mbox{\rm Dim}^*\mu $.

\subsection{Sums, products, convolutions, projections of measures}
What we present in this subsection was in the first unpublished version of
\cite{Fan1994}. Some part was restated by in  \cite{Tamashiro1995} and some part
was used in \cite{Fan1994b,Fan1994c}.

{\em Sum.}
The absolute continuity $\nu\ll \mu$ is a partial order on the space of positive Borel
measures $M^+(X)$ on the space $X$. 
Using the continuity-singularity criterion, it is easy to see that
\begin{equation}\label{dim-equiv}
    0<\nu \ll \mu \Rightarrow \dim_*\mu \le \dim_*\nu \le \dim^* \nu \le \dim^* \mu.
\end{equation}
The relation $\nu\sim \mu$ (meaning $\nu\ll \mu \ll \nu$) is an equivalent relation.
Since two equivalent measures have the same lower and upper Hausdorff dimensions,
both $\dim_*\mu$ and $\dim^*$ are well defined for equivalent classes.

Given a family of positive measures $\{\mu_i\}_{i\in I}$ which is bounded under the order
$\ll$, we denote  its supremum by $\bigvee_{i\in I} \mu_i$. If the family
is finite, we have $\bigvee_{i\in I} \mu_i \sim \sum_{i\in I} \mu_i$. Such equivalence
also holds when the family is countable. In general, there is a countable sub-family $\mu_{i_k}$
such that
$
     \bigvee_{i\in I} \mu_i \sim \sum_k \mu_{i_k}.
$
This is what we mean by sum of measures.

\begin{theorem} If $\{\mu_i\}_{i\in I}$ is a bounded family of measures in $M^+(X)$, we have
$$
    \dim_* \bigvee_{i\in I} \mu_i = \inf_{i \in I} \dim_*\mu_i,\quad
    \dim^* \bigvee_{i\in I} \mu_i = \sup_{i \in I} \dim^*\mu_i.
$$
\end{theorem}

Let us see how to prove the first formula for a family of two measures $\mu$ and $\nu$.
For any $\alpha < \dim_* (\mu+\nu)$, we have
$
     U_\alpha^\mu(x) + U_\alpha^\mu(x) = U_\alpha^{\mu+\nu}(x)<\infty \quad (\mu+\nu)\!-\!a.e.
$
This implies $\dim_*(\mu+\nu) \le \min\{\dim_*\mu, \dim_* \nu\}$.
The inverse inequality follows from the fact that if
$ \beta < \min\{\dim_*\mu, \dim_* \nu\}$, then
$
    \mu \ll H^\beta \ \mbox{\rm and}\ \nu\ll H^\beta
$
which implies $\mu+\nu \ll H^\beta$.
\medskip

Let us consider now the infimum $\bigwedge_{i \in I}$ of a family of
measure $\{\mu_i\}_{i\in I}$.
Recall that by definition we have
$$
  \bigwedge_{i\in I} \mu_i \sim \bigvee_{\mu: \forall i\in I, \mu\ll \mu_i} \mu.
$$
For a family of two measures $\mu$ and $\nu$,
we have
$
      \mu\bigwedge\nu \sim \frac{d\mu}{d \nu} \bigvee \frac{d\nu}{d\mu}.
$

\begin{theorem} If $\{\mu_i\}_{i\in I}$ is  family of measures in $M^+(X)$ such that
$\bigwedge_{i \in I} \mu_i\not=0$, we have
$$
    \sup_{i\in I}\dim_*  \mu_i
    \le \dim_* \bigwedge_{i \in I} \mu_i
    \le \dim^* \bigwedge_{i \in I} \mu_i
    \le  \inf_{i\in I}\dim^*  \mu_i.
    $$
\end{theorem}

Consequently, $\bigvee_{i \in I} \mu_i$ (resp. $\bigwedge_{i \in I} \mu_i$)
is unidimensional if and only if all measures $\mu_i$ are unidimensional and have the same dimension.

\medskip
{\em Product.}
Let $(X,\delta_X)$ and $(Y, \delta_Y)$ be two Polish spaces. Then
$d:=\delta_X \vee \delta_Y$ is a compatible metric on the product space.
Let $\mu \in M^+(X)$ and $\nu\in M^+(Y)$. Let us consider the product measure
$\mu \otimes\nu$.

\begin{theorem} For $\mu \in M^+(X)$ and $\nu\in M^+(Y)$, we have
$$
   \dim_* \mu\otimes \nu \ge \dim_*\mu + \dim_*\nu, \qquad
   \dim^* \mu\otimes\nu \ge \dim^*\mu + \dim^*\nu.
$$
\end{theorem}

We say a measure $\mu\in M^+(X)$ is regular  if $\underline{D}(\mu, x)= \overline{D}(\mu, x)$
$\mu$-a.e.

\begin{theorem} If $\mu \in M^+(X)$ or $\nu\in M^+(Y)$ is regular, we have
$$
   \dim_* \mu\otimes \nu = \dim_*\mu + \dim_*\nu, \qquad
   \dim^* \mu\otimes\nu = \dim^*\mu + \dim^*\nu.
$$
\end{theorem}
\medskip

{\em Convolution.}  Assume that the Polish space $X$ is a locally compact abelian group $G$.
Assume further that $G$ satisfies the following hypothesis
$$
     I_\alpha^{\mu*\nu} \le C(\alpha, \nu) I_\alpha^\mu
$$
for all measures $\mu, \nu \in M^+(X)$, where $C(\alpha, \nu)$ is a constant independent of $\mu$.
For example, if $G=\mathbb{R}^d$, we have
$$
     I_\alpha^\mu = \int |\xi|^{-(d-\alpha)} |\widehat{\mu}(\xi)|^2 d\xi
$$
which implies that the hypothesis is satisfied by $\mathbb{R}^d$. The hypothesis is also
satisfied by the group $\prod_{n=1}^\infty \mathbb{Z}/m_n\mathbb{Z}$ \cite{Fan1989c}.

\begin{theorem} For any measures $\mu, \nu \in M^+(G)$, we have
$$
   \dim_* \mu*\nu \ge \max\{\dim_*\mu, \dim_*\nu\}, \qquad
   \dim^* \mu*\nu \ge \max\{\dim^*\mu, \dim^*\nu\}.
$$
\end{theorem}

For a given measure $\mu\in M^+(G)$, we consider the following two subgroups
$$
    H_-:=\{t \in G: \mu \ll \mu *\delta_t\},
    \qquad  H_+:=\{t \in G: \mu*\delta_t \ll \mu\}.
$$
We could call $H = H_-\cap H_+$ the quasi-invariance group of $\mu$.

\begin{theorem} Under the above assumption,  if  $\nu(H_-)>0$, we have the equality
$
   \dim_* \mu*\nu = \dim_*\mu
$ and consequently $\dim_*\nu\le \dim_*\mu$ and $\dim_H H_-\le \dim_*\mu$; if  $\nu(H_+)>0$, we have the equality
$
   \dim^* \mu*\nu = \dim^*\mu
$ and consequently $\dim^*\nu\le \dim^*\mu$ and $\dim_H H_+\le \dim^*\mu$.
\end{theorem}

\medskip
{\em Projection.}
Let $\mu \in M^+(\mathbb{R}^2)$ with $\mathbb{R}^2=\mathbb{C}$. Let $L_\theta$ ($0\le \theta<2\pi$) be the line
passing the origin and having angle $\theta$ with the line of abscissa. The orthogonal projection
$P_\theta$
on $L_\theta$ is defined by $P_\theta(x)= \langle x, e^{i\theta} \rangle$.
Let $\mu_\theta: = \mu \circ P_\theta^{-1}$ be the projection of $\mu$ on $L_\theta$.

\begin{theorem} Let $\mu \in M^+(\mathbb{R}^2)$. For almost all $\theta$,
we have $\dim_*\mu_\theta = \dim_* \mu \wedge 1$ and $\dim^*\mu_\theta = \dim^* \mu \wedge 1$.
\end{theorem}

\subsection{Ergodicity and dimension}
L. S. Young \cite{Young1982} considered diffeomorphisms of surfaces leaving invariant an ergodic Borel probability measure $\mu$. She proved that $\mu$ is exact and found  a formula relating $\dim \mu$ to the entropy and Lyapunov exponents of $\mu$. One of the main problems in the interface of dimension theory and dynamical systems is the Eckmann-Ruelle conjecture on the dimension of hyperbolic ergodic measures: the local dimension
of every hyperbolic measure invariant under a $C^{1+\alpha}$-diffeomorphism exists almost everywhere.
This conjecture was proved by L. Barreira,  Y. Pesin and J. Schmeling~\cite{BPS} based on the fundamental
fact that such a measure possesses asymptotically "almost" local product structure.
But, in general, the ergodicity of the measure doesn't imply that the measure is exact~\cite{Cutler}.
In \cite{Fan1994b}, $D$-ergodicity and unidimensionality were studied.

\section{Oriented walks and Riesz products}\label{Sect_OW}
\subsection{Oriented walks}
Let $(\epsilon_n)_{n\ge 1} \subset [0, 2\pi)^\mathbb{N}$ be a sequence of angles. For
$n\ge 1$, define
$$
     S_n(\epsilon) = \sum_{k=1}^n e^{i(\epsilon_1 + \epsilon_2 + \cdots + \epsilon_k)}.
$$
We call $(S_n(\epsilon))_{n\ge 1}$ an oriented walk on the plan $\mathbb{C}$.
In his book \cite{Feller} (vol. 1, pages 240-241), Feller mentioned  a
model describing the length of long polymer molecules in chemistry.
It is a random chain consisting of $n$ links, each of unit length,  and the angle between two consecutive links is $\pm \alpha$ where $\alpha$ is a positive constant.
Then the distance $L_n$
from the beginning to the end of the chain can be expressed by
$$
    L_n = |S_n(\epsilon)|
$$
where $(\epsilon_n)$ is an i.i.d. sequence of random variables taking values in
$\{-\alpha, \alpha\}$.
If $\alpha=0$,
$L_n=n$ is deterministic. If $0<\alpha<2\pi$, the random variable $L_n$ is not expressed as sums of
independent variables. However
Feller succeeded in computing the second order moment of $L_n$.
It is actually proved in \cite{Feller} that $\|L_n\|_2$ is of order
$\sqrt{n}$. More precisely, for $0<\alpha<2\pi$ we have
$$
    \mathbb{E} L_n(\alpha)^2 = n \frac{1 + \cos \alpha}{1-\cos \alpha }
     - 2 \cos \alpha  \frac{1 - \cos^n \alpha}{(1-\cos \alpha)^2 }.
$$
Observe that
$$
\mathbb{E}L_n^2 = \frac{1-(-1)^n}{2} \ \mbox{\rm if}\ \alpha=\pi;
\qquad \mathbb{E}L_n^2\sim n \frac{1+\cos \alpha}{1- \cos \alpha}
\ \mbox{\rm if}\ 0<\alpha<2\pi, \alpha\not=\pi.
$$

What is the behavior of $S_n(\epsilon)$ as $n\to \infty$ for individuals $\epsilon$ ?
We could study the behavior from the multifractal point of view.
Let us consider a more general setting.
Fix $d\in \N^*$. Let $\tau \in GL(\R^d)$, $v \in \mathbb{R}^d$ and $A$ a finite subset of
$\Z$. For any $x=(x_n) \in \mathbb{D}:=A^\N$,
we define the oriented walk
$$
S_0(x)=v, \qquad S_n(x)=\sum_{k=1}^n \tau^{x_1+x_2+\cdots +x_k}v.
$$
For $\alpha\in \R^d$, we define the $\alpha$-{\em level set}
$$E_{\tau}(\alpha):=\left\{x\in  \mathbb{D} : \lim_{n\to\infty }\frac1n S_n(x)=\alpha\right\}.$$
Let $L_\tau:=\{\alpha\in\R^d : E_{\tau}(\alpha)\neq \emptyset \}$.

The following two cases were first studied in   \cite{Fan2000}.
Case 1: $d=2$, $\tau = -1$ and $A=\{0, 1\}$; Case 2: $d=1$, $\tau = e^{i \pi/2}$ and $A=\{-1, 1\}$.

\begin{theorem}\cite{Fan2000}.
In the first case, we have  $L_\tau=[-1,1]$ and for $\alpha\in L_\tau$ we have
$$\dim_H E_{\tau}(\alpha)= \dim_P E_{\tau}(\alpha)=H\left(\frac{1+\alpha}{2}\right).$$
In the seconde case, we have  $L_\tau=\{z=a+ib : |a|\leq 1/2, |b|\leq 1/2 \}$ and for $\alpha=a+bi\in L_\tau$ we have
$$
\dim_H E_{\tau}(\alpha)= \dim_P E_{\tau}(\alpha)=\frac{1}{2\log 2}\left[H\left(\frac{1}{2}+a\right)+H\left(\frac{1}{2}+b\right)\right].$$
\end{theorem}

This theorem was proved by using Riesz products which will be described in the following subsection.

A new construction of measures allows us to deal with a class of oriented walks.
We assume that  $\tau \in GL(\R^d)$ is idempotent. That is to say
$\tau^p = Id$ for some integer $p>1$ (the case $p=1$ is trivial). The least $p$ is called the order
 of $\tau$. The above two cases are special
cases. In fact,  $\theta= -1$ is idempotent with order $p=2$ and
 $\theta=e^{i\pi/2}$ is idempotent with order $p=4$. The following rotations
 in $\mathbb{R}^3$
$$
\tau_1 =
\left( \begin{array}{ccc}
0 & 0 & 1 \\
1 & 0 & 0 \\
0 & 1 & 0
\end{array} \right), \qquad
\tau_2 =
\left( \begin{array}{ccc}
0 & 0 & -1 \\
1& 0 & 0 \\
0 & 1 & 0
\end{array} \right)
$$
are idempotent with order respectively equal to $3$ and $6$. Also remark that
$\tau \in SO_3(\mathbb{R})$ and $\tau_2 \in O_3(\mathbb{R})\setminus SO_3(\mathbb{R})$.

Since $\tau^p=Id$, the sum $x_1 + x_2+\cdots +x_k$ in the definition of $S_n(x)$
can be made modulo $p$.
For $s\in\R^d$, we define a $p\times p$-matrix $M_s = (M_s(i, j))$: for $(i, j) \in \mathbb{Z}/p\mathbb{Z}
\times \mathbb{Z}/p\mathbb{Z}$ define
 $$
 M_s(i,j) = 1_A(j-i) \exp[\langle s, \tau^j\rangle v].
 $$
where $\langle \cdot, \cdot \rangle$ denotes the scalar product in $\R^d$. It is clear that
$M_s$ is irreducible iff $M_0$ is so. We consider $A$ as a subset (modulo $p$) of $\mathbb{Z}/p\mathbb{Z}$. It is easy to see that $M_0$ is irreducible iff $A$
generates the group $\mathbb{Z}/p\mathbb{Z}$.

Assume that $A$ generates the group $\mathbb{Z}/p\mathbb{Z}$.
Then $M_s$ is irreducible and by the Perron-Frobenius theorem, the spectral radius
$\lambda(s)$ of $M_s$ is a simple eigenvalue and there is a unique corresponding
probability eigenvector $t(s) = (t_s(0), t_s(1) \cdots, t_s(p-1)$.
Let
$$
  P(s) = \log \lambda(s).
$$
It is real analytic and strictly convex function on $\R^d$. We call it the pressure function
associated to the oriented walk.

\begin{theorem} \label{Thm_OW}\cite{FW}
Assume $\tau$ is idempotent with order $p$ and $A$ generates the group
$\mathbb{Z}/p \mathbb{Z}$. Then $L_\tau = \overline{\{\triangledown P(s) : s\in \R^d \}}$
and for  $\alpha\in \bigtriangleup$, we have
$$\dim_H E_{\tau}(\alpha)= \dim_P E_{\tau}(\alpha)=\frac{1}{\log p}\inf_{s\in\R^d}\{P(s)-\langle s,\alpha\rangle\}=
\frac{P(s_\alpha)-\langle s_\alpha,\alpha\rangle}{\log p},$$
where $s_\alpha$ is the unique $s\in\R^d$ such that $\triangledown P(s)=\alpha$.
\end{theorem}

\subsection{Riesz products}

%
Theorem~\ref{Thm_FLM} was proved by using Riesz products.
While Hausdorff introduced the Hausdorff dimension (1919), F. Riesz constructed a class
of continuous but singular measures on the circle (1918), called Riesz products. Riesz products are used as tool in harmonic analysis and some of them are Gibbs measures in the sense of dynamical systems.

Let us recall the definition of Riesz product on a compact abelian group $G$, due to
Hewitt-Zuckerman \cite{HZ} (1966).
Let $\widehat{G}$ be the dual group of $G$. 
A sequence of characters $\Lambda = (\gamma_n )_{n
  \geq 1} \subset \widehat{G}$ is said to be {\em dissociated} if for any $n
\ge 1$, the following characters are all distinct:
$$
\gamma_1^{\epsilon_1} \gamma_2^{\epsilon_2}\cdots
\gamma_n^{\epsilon_n}
$$
where $\epsilon_j \in
 \{-1, 0, 1\}$ if $\gamma_j$ is not of order $2$,
 or $\epsilon_j \in
\{0,1\}$ otherwise.
Given such a dissociated sequence $\Lambda = (\gamma_n )_{n \geq 1}$ and a sequence of complex numbers $a = (a_n)_{n \geq 1}$ such that $|a_n | \leq 1$, we can define a probability measure on
 $G$, called {\em Riesz product},
\begin{equation}\label{Riesz}
  \mu_a =\prod_{n=1}^\infty \bigl(1 + \mbox{\rm Re}\  a_n \gamma_n (t)\bigr)
\end{equation}
as the weak* limit of  $\prod_{n=1}^N \bigl(1 + \mbox{\rm Re}\ a_n \gamma_n
(t)\bigr) \text{d}t$ where $\text{d}t$ denotes the Haar measure on $G$.

A very useful fact is that the Fourier coefficients of the Riesz product
 $\mu_a$ can be explicitly expressed in term of the
coefficients $a_n$'s:
$$
\widehat{\mu}_a (\gamma) =
\prod_{k=1}^n a_k^{(\epsilon_k)} \ \mathrm{if}  \  \gamma=\gamma_1^{\epsilon_1} \gamma_2^{\epsilon_2}\cdots
\gamma_n^{\epsilon_n}, \quad \widehat{\mu}_a (\gamma) =0 \ \mathrm{otherwise}.
$$
where $a_n^{(\epsilon)} = 1,
a_n/2$ or $\bar{a}_n/2$ according to $\epsilon = 0, 1$ or $-1$. Consequently
the sequence $\{ \gamma_n - \bar{a}_n/2 \}_{n\ge 1}$
is an orthogonal system in $L^2(\mu_a)$. Here are some properties of $\mu_a$.

\begin{theorem}[\cite{Zygmund1968}] \label{ACS} The measure $\mu_a$ is
  either absolutely continuous or singular (with respect to the Haar measure) according to $\sum_{n=1}^\infty |a_n|^2<\infty$ or
  $=\infty$.
\end{theorem}

\begin{theorem} [\cite{Fan1993},\cite{Peyriere1990}] \label{FP}
  Let $\{\alpha_n\}$ be a sequence of complex numbers. The orthogonal series
 $\sum_{n=1}^\infty \alpha_n (\gamma_n(t) - \bar{a}_n/2)$
  converges $\mu_a$-everywhere iff
  $\sum_{n=1}^\infty |\alpha_n|^2 <\infty.$
\end{theorem}

The proof of Theorem~\ref{FP} in \cite{Fan1993} involved the following Riesz product
with $\omega =(\omega_n) \in
G^\mathbb{N}$ as phase translation:
\begin{equation}\label{RandomRiesz}
  \mu_{a,\omega} =\prod_{n=1}^\infty (1 + \mbox{\rm Re} \ a_n \gamma_n
  (t+\omega_n)).
\end{equation}
Actually $\mu_{a, \omega}$ was considered as a random measure and $(\omega_n)$
was considered as an i.i.d. random sequence with Haar measure as common
probability law.

When two Riesz products $\mu_a$ and $\mu_b$ are singular or mutually absolutely
continuous ?
It is a unsolved problem. Bernoulli infinite product measures can be viewed
as Riesz products on the group $(\mathbb{Z}/m\mathbb{Z})^{\mathbb{N}}$. For these
Bernoulli infinite measures, the Kakutani dichotomy theorem \cite{Kakutani} applies and there is a complete
solution. But there is no complete solution for other groups.

The classical Riesz products are of the form
\begin{equation}\label{RealRiesz}
   \mu_a = \prod_{k=1}^\infty (1 + \mbox{\rm Re} \ a_n e^{i \lambda_k x})
\end{equation}
where $\{\lambda_n\}\subset \mathbb{N}$ is a lacunary sequence in the sense that
$\lambda_{n+1} \ge 3 \lambda_n$. J. Peyri\`ere has first studied the
lower and upper dimensions of $\mu_a$, without introducing the notion of dimension
of measures. 
Let us mention an estimation
for the energy integrals of $\mu_a$ \cite{Fan1989c}:
$$
   \int\int \frac{d \mu_a(x)d\mu_a(y)}{|x-y|^\alpha}
   \approx \lambda_1^{\alpha -1}|a_1|^2
   + \sum_{n=2}^\infty \lambda_n^{\alpha -1}|a_n|^2 \prod_{k=1}^{n-1} \left(1 + \frac{|a_n|^2}{2}\right).
$$

\subsection{Evolution measures}
The key for the  proof of Theorem~\ref{Thm_OW} is the construction of the following
measures on $A^\mathbb{N}$, which describe the evolution of the oriented walk. It is similar to Markov measure
but it is not. It plays the role of Gibbs measure but it is not Gibbs measure
either.

Recall that $M_s t(s) = \lambda(s) t(s)$. In other words, for every $i \in
\mathbb{Z}/p\Z$ we have
$$
\lambda(s) t_i(s) =\sum_{j} 1_A(j -i) t_j(s) \exp[\langle s, \tau^j v\rangle].
$$
Denote, for $a\in A$ and for $(x_1, \cdots, x_{k+1}) \in A^{k+1}$,
$$
\pi(a)=\frac{t_a(s)}{\sum_{b \in A}t_{b}(s)};
$$
$$
Q_k(x_1,x_2,\cdots, x_{k+1})=
\frac{t_{x_1 + x_2+\cdots +x_{k+1}}(s) \exp[\langle s,\tau^{x_1 + x_2+\cdots +x_{k+1}}v \rangle]} {\lambda(s) t_{x_1 + x_2+\cdots +x_{k}}(s)}.$$
Then
we define a probability measure $\mu_s$ on $A^\mathbb{N}$ as follows. For any word $x_1x_2\cdots x_n\in A^n$, let
$$\mu_s([x_1x_2\cdots x_n])=\pi(x_1)Q_1(x_1,x_2)Q_2(x_1,x_2,x_3)\cdots Q_{n-1}(x_1,x_2,\cdots, x_n).$$

For $x = (x_n) \in A^\mathbb{N}$, let
$$
     w_n(x) = x_1 + x_2 + \cdots + x_n      \quad (\!\!\!\!\!\mod p).
$$
The mass $\mu_s([x_1x_2\cdots x_n])$ and the partial sum $S_n(x)$ are directly related as follows.
$$\log\mu_s([x_1x_2\cdots x_n])=\left\langle s,S_n(x)\right\rangle-(n-1)\log \lambda(s)-\log \sum_{a\in A} t_{a}(s)+\log t_{w_n(x)}(s).$$
As the $t_i(s)$'s are bounded, we deduce that the following relation between the measure $\mu_s$
and the oriented walk $S_n$.

\begin{proposition}\label{prop local dimension1}
For any $x\in \mathbb{D}$, we have
$$\log\mu_s([x_1x_2\cdots x_n]) - \left\langle s,S_n(x)\right\rangle=
- n \log \lambda(s)+ O(1).$$
\end{proposition}




\section{Multiple Birkhoff averages}\label{Sect_MBA}

 Let $\mathcal{A}=\{0, 1,2, \cdots, m-1\}$
be a set of $m$ symbols ($m\ge 2$). Denote $\Sigma_m=\mathcal{A}^{\mathbb{N}^*}$. Let $q\ge 2$ be a integer.
 Fan, Schmeling and Wu made a forward step in \cite{FSW} by obtaining a Hausdorff spectrum
of multiple ergodic averages for a class of potentials.
They consider  an arbitrary function $\varphi: \mathcal{A}^d \to \mathbb{R}$  and study the sets
$$E(\alpha)=\left\{x\in \Sigma_m : \lim_{n\to\infty}A_n\varphi(x)=\alpha\right\}$$
for $\alpha\in \R$, where
\begin{equation}\label{GMEA}
A_n\varphi(x)=\frac{1}{n}\sum_{k=1}^n
  \varphi(x_k, x_{qk}, \cdots, x_{q^{d-1}}).
  \end{equation}
Let
$$\alpha_{\min}=\min_{a_1,\cdots,a_{d}\in
\mathcal{A} }\varphi(a_1,\cdots,a_{d}), \quad
\alpha_{\max}=\max_{a_1,\cdots,a_{d}\in
\mathcal{A} }\varphi(a_1,\cdots,a_{d}).$$ It is assumed that
$\alpha_{\min}<\alpha_{\max}$ (otherwise $\varphi$ is constant and
the problem is trivial). A key ingredient of the proof is a class of measures constructed by Kenyon, Peres and Solomyak \cite{KPS} that we call telescopic product measures. In \cite{FSW}, a nonlinear thermodynamic
formalism was developed.

\subsection{
Thermodynamic formalism}
\label{M-sets}\ \\

The Hausdorff dimension of $E(\alpha)$ is determined through the following thermodynamic formalism.
Let $\mathcal{F}(\mathcal{A}^{d-1}, \mathbb{R}^+)$ be the cone of
functions defined on $\mathcal{A}^{d-1}$ taking non-negative real values. For any $s \in
\mathbb{R}$, consider the  transfer operator $\mathcal{L}_s$ defined on
$\mathcal{F}(\mathcal{A}^{d-1}, \mathbb{R}^+)$ by
\begin{equation}\label{transer-operator}
\mathcal{L}_s \psi (a)
= \sum_{j \in \mathcal{A}} e^{s \varphi(a, j)}
\psi (Ta, j)
\end{equation}
where $T: \ \mathcal{A}^{d-1}\to \mathcal{A}^{d-2}$ is defined by $T(a_1,\cdots,a_{d-1})=(a_2,\cdots,a_{d-1})$.
Then define the non-linear operator $\mathcal{N}_s$ on $\mathcal{F}(\mathcal{A}^{d-1}, \mathbb{R}^+)$ by
$
     \mathcal{N}_s \psi (a)= (\mathcal{L}_s \psi (a))^{1/q}.
$
It is
proved in \cite{FSW} that the equation
\begin{equation}\label{transer_equation}
     \mathcal{N}_s \psi_s = \psi_s
\end{equation}
admits a unique strictly positive solution $\psi_s=\psi_s^{(d -1)} : \mathcal{A}^{d-1}\to
\mathbb{R}_+^*$.
Extend the function $\psi_s$  onto $\mathcal{A}^{k}$ for all $1\le k \le d -2$ by induction:
\begin{equation}\label{transer_equation2}
\psi_s^{(k)} (a)=\left(\sum_{j \in \mathcal{A}} \psi_s^{(k+1)} (a, j)\right)^{\frac{1}{q}}, \ \ (a\in \mathcal{A}^{k}).
\end{equation}
For simplicity, we  write $\psi_s(a)=\psi_s^{(k)}(a)$ for $a\in \mathcal{A}^k$ with $1\leq k\leq d-1$.
Then the pressure function is defined by
\begin{equation}\label{pressure function}
P_{\varphi}(s) = (q-1)q^{d-2} \log \sum_{j\in \mathcal{A}}\psi_s(j).
\end{equation}
It is  proved \cite{FSW} that $P_{\varphi}(s)$ is an analytic and convex function of $s\in
\mathbb{R}$ and even strictly convex since $\alpha_{\min}
<\alpha_{\max}$. The Legendre transform of $P_{\varphi}$ is defined
as $$ P^*_{\varphi}(\alpha)=\inf_{s\in \mathbb{R}}(P_{\varphi}(s)-s\alpha).
$$
We denote by $L_{\varphi}$ the set of levels $\alpha\in \R$ such that
$E(\alpha)\neq \emptyset$.
    \medskip

\begin{theorem} {\rm (\cite{FSW})}\label{thm principal}
We have
$L_{\varphi}=[P'_{\varphi}(-\infty),P'_{\varphi}(+\infty)].$
 If
$\alpha =P'_{\varphi}(s_\alpha)$ for some $s_\alpha \in
\R\cup\{-\infty,+\infty\}$, then $E(\alpha)\neq \emptyset$ and the
Hausdorff dimension of $E(\alpha)$ is equal to
$$\dim_H E(\alpha)=\frac{P_{\varphi}^*(\alpha)}{q^{d-1}\log m}.$$
\end{theorem}

\mk
Similar results hold for vector valued functions $\varphi$ \cite{FSW}.
Y. Peres and B. Solomyak \cite{PS2} have obtained a result for the special case
$\varphi(x,y) = x_1y_1$ on $\Sigma_2$.
Y. Kifer \cite{Kifer} has obtained a result on the multiple recurrence sets
for some frequency of product form.

Let us consider two examples.
Let $q=2$ and $\ell=2$ and  let $\varphi_1(x,y)=x_1y_1$ (See Figure 1) and $\varphi_2(x,y) = (2x_1-1)(2y_1-1)$ (see Figure 2)
be two potentials on $\Sigma_2$. The invariant spectra (see \S \ref{InvSpectrum}) are also shown in the figures.

\begin{figure}
\centering
\includegraphics[width=6.8cm]{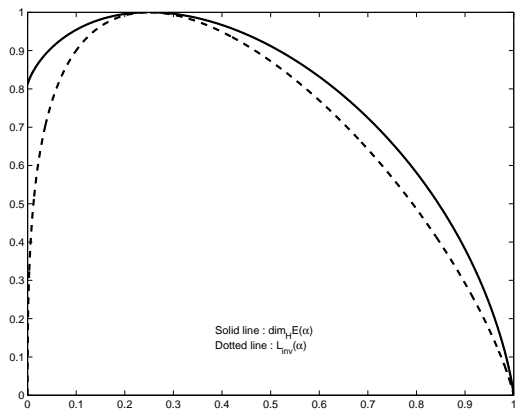}
\caption{Spectra $\alpha\mapsto \dim_HE(\alpha)$ and $\alpha\mapsto F_{\rm inv}(\alpha)$ for $\varphi_1$.}
\label{figure 1}
\end{figure}

\begin{figure}
\centering
\includegraphics[width=6.8cm]{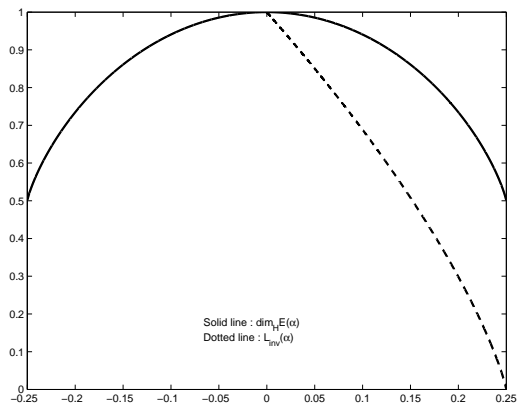}
\caption{Spectra $\alpha\mapsto \dim_HE(\alpha)$ and $\alpha\mapsto F_{\rm inv}(\alpha)$ for $\varphi_2$.}
\label{figure 2}
\end{figure}


\subsection{Telescopic product measures}
One of the key points in the proof of the Hausdorff spectrum
 (Theorem \ref{thm principal}) is the observation that the coordinates $x_1,\cdots , x_n,\cdots $ of $x$ appearing in the definition of ${A}_n\varphi(x)$ share the following independence.  Consider the partition of $\N^*$:
 $$
 \N^*=\bigsqcup_{i\geq 1,q\nmid i}\Lambda_i\ \ {\rm with}\ \Lambda_i=\{iq^j\}_{j\ge 0}.
 $$
 Observe that if $k=iq^j$ with $q\nmid i$, then $\varphi(x_k,x_{kq},\cdots ,x_{kq^{d-1}})$ depends only on $x_{|_{\Lambda_i}}$, the restriction of $x$ on $\Lambda_i$. So the summands in the definition of ${A}_n\varphi(x)$ can be put into different groups, each of which depends on one restriction $x_{|_{\Lambda_i}}$. For this reason, we decompose $\Sigma_m$ as follows:
 $$
 \Sigma_m =\prod_{i\geq1,q\nmid i}\mathcal{A}^{\Lambda_i}.
 $$

Telescopic product measures are now constructed as follows.
 Let $\mu$ be a probability measure on $\Sigma_m$. Notice that $\mathcal{A}^{\Lambda_i}$ is nothing but a copy of $\Sigma_m$.
 We consider $\mu$ as a measure on $\mathcal{A}^{\Lambda_i}$ for every $i$ with $q\nmid i$. Then we define the infinite product measure $\P_\mu$ on   $\prod_{i\geq1,q\nmid i}\mathcal{A}^{\Lambda_i}$ of the copies of $\mu$. More precisely,  for any word $u$ of length $n$ we define
 $$
 \P_{\mu}([u])=\prod_{i\leq n,q\nmid i}\mu([u_{|_{\Lambda_i}}]),
 $$
 where $[u]$ denotes the cylinder of all sequences starting with $u$. The probability measure  $\P_{\mu}$ is  called  {\em telescopic product measure}.  Kenyon, Peres and Solomyak \cite{KPS} have first constructed  these measures.
 \medskip

 The Hausdorff dimension of every telescopic product measure is computable.

 \begin{theorem}{\rm (\cite{KPS, FSW})}\label{dim-meas} For any given measure $\mu$,
the telescopic product measure $\mathbb{P}_\mu$ is exact and its dimension is equal to
  $$\dim_{H}\P_{\mu}=\dim_P \P_{\mu}=\frac{(q-1)^2}{\log m}\sum_{k=1}^{\infty}\frac{H_k(\mu)}{q^{k+1}}$$
  where
  $H_k(\mu)=-\sum_{a_1,\cdots,a_k\in S}\mu([a_1\cdots a_k])\log \mu([a_1\cdots a_k]).
  $
\end{theorem}

\subsection{Dimension formula of Ruelle-type}

The function $\psi_s$ defined by (\ref{transer_equation}) and (\ref{transer_equation2})
determine a special telescopic product measure which plays the role of Gibbs measure
in the proof of the Hausdorff spectrum.

 First we define  a $(d-1)$-step Markov measure on $\Sigma_m$, which will
  be denoted by  $\mu_s$,  with the initial law
\begin{equation}\label{def measure 1}
\pi_s([a_1,\cdots,a_{d-1}])=\prod_{j=1}^{d-1}\frac{\psi_s(a_1,\cdots,a_j)}{\psi_s^q(a_1,\cdots,a_{j-1})}
\end{equation}
and the transition probability
\begin{equation}\label{def measure 2}
Q_s\left([a_1,\cdots,a_{d-1}],[a_2,\cdots,a_{d}]\right)=e^{s\varphi(a_{1},\cdots,a_{d})} \frac{\psi_s(a_2,\cdots,a_{d})}{\psi_s^q(a_{1},\cdots,a_{d-1})}.
\end{equation}
The corresponding telescopic product measure $\mathbb{P}_{\mu_s}$
is proved to be a dimension maximizing measure of $E(\alpha)$ if $s$ is chosen to be the solution of
$
    P'_\varphi(s) =\alpha.
$
The dimension of $\mathbb{P}_{\mu_s}$ is simply expressed by the pressure function. In other words, we have the following formula of Ruelle-type.
\begin{theorem}{\rm (\cite{FSW})} \label{prop formula dim P_s}
For any $s\in \R$, we have
$$\dim_H\P_{\mu_s}=\frac{1}{q^{d-1}}[P_\varphi(s)-sP'_\varphi(s)].$$
\end{theorem}

\subsection{
Multiplicatively invariant sets}
\label{M-sets}\ \\

Kenyon, Peres and Solomyak \cite{KPS} were able to compute both the Hausdorff dimension and the box dimension of $X_2$, already considered in \S \ref{X_2},  and of a class of generalizations of $X_2$.
Peres, Schmeling, Seuret and Solomyak \cite{PSSS} generalized the results to a more general class of sets.

Recall that $\Sigma_m=\{0,1,..., m-1\}^{\N^*}$,
$q$ is an integer greater than 2 and $\Lambda_i : =\{i\,q^n: n\geq 0\}$.
For
any subset $\Omega \subset \Sigma_m$, define
\begin{equation}
\label{defXOmega}
X_\Omega := \{x=(x_k)_{k\geq 1}\in \Sigma_m: \forall \, i\neq 0 \!\! \!\mod q, \ x_{|\Lambda_i} \in \Omega\}.
\end{equation}
We get $X_\Omega =X_2$ when
$q=2$ and $\Omega$ is the Fibonacci set $\{y\in \Sigma_m: \forall \, k\geq 1, \ y_ky_{k+1}=0\}$.

The set $X_\Omega$ is not shift-invariant but is multiplicatively invariant in the sense that
$M_r X_\Omega \subset X_\Omega$ for every integer $r\in \mathbb{N}^*$ where
$M_r$ maps $(x_n)$ to $(x_{rn})$.

The generating set $\Omega$ has a tree of prefixes, which is a directed graph $\Gamma$. The  set  $V(\Gamma)$
 of vertices consists of all possible prefixes of finite length in $\Omega$, i.e.
$$V(\Gamma) := \bigcup_{k\geq 0} \mbox{Pref}_k(\Omega) ,$$
where $\mbox{Pref}_0(\Omega) =\{\emptyset\}$,
  $\mbox{Pref}_k(\Omega) := \{u\in \{0,1,..., m-1\}^k: \Omega \cap [u]\neq \emptyset\}$.
There is a directed edge   from a vertex $u$ to another $v$ if and only if  $v=ui$ for some $i\in \{0,1,..., m-1\}$.

\begin{theorem}\label{thKPS}{\rm (\cite{KPS})}
There exists a unique vector  $\bar t=(t_v)_{v\in \Gamma} \in [1,m^{\frac{1}{q-1}}]^{V(\Gamma)}$
 defined on the tree such that
\begin{equation}
\label{deftoptimal}
 \forall v\in V(\Gamma),  \ \ \ (t_v)^q = \sum_{i\in\{0,1,...,m-1\}: \ vi\in \Omega} \ t_{vi}.
 \end{equation}
The Hausdorff dimension and the box dimension of $X_\Omega$ are respectively equal to
\begin{eqnarray}
\dim_H (X_\Omega )  &  =  & (q-1) \log_m t_\emptyset\\
\dim_B (X_\Omega )  &  =  & (q-1)^2 \sum_{k=1}^{+\infty} \frac{ \log_m  \left| \text{{\em Pref}}_k(\Omega)  \right|}{q^{k+1}}.
\end{eqnarray}
The two dimensions coincide if and only if the tree $\Gamma$ is spherically symmetric, i.e. all prefixes of length $k$ in $\Omega$ have the same number of continuations of length $k+1$ in $\Omega$.
\end{theorem}

 The vector $\bar{t}$ defines a measure $\mu$ on $\Omega\subset \{0, 1, \cdots, m-1\}^{\mathbb{N}^*}$.
 Then a telescopic product measure can be  built on $X_\Omega$.
It is proved in \cite{KPS} that there is a maximizing measure on $X_\Omega$ of this form.

\mk

A typical example of the class of sets studied by Peres, Schmeling, Solomyak and Seuret \cite{PSSS} is
$$X_{2,3}=\{x\in \Sigma_m: x_kx_{2k}x_{3k} = 0\}.$$
The construction of the sets is as follows. Let   $\kappa \ge 1$ be an integer and let $p_1,\ldots,p_\kappa$ be $\kappa$ primes, which generates a semigroup $S$ of $\mathbb{N}^*$:
$$S=\langle p_1, p_2, \cdots, p_\kappa\rangle = \{p_1^{\alpha_1}p_2^{\alpha_2} \cdots p_\kappa^{\alpha_\kappa}: \alpha_1, ..., \alpha_\kappa\in \N\}.$$
 The elements of $S$ are arranged in increasing order and denoted by $\ell_k$ the $k$-th element of
$
S = \{\ell_k\}_{k=1}^\infty: 1=\ell_1 < \ell_2 < \cdots.
$
Define
\begin{equation} \label{eq-fs}
\gam(S):= \sum_{k=1}^\infty \frac{1}{\ell_k}\,.
\end{equation}
Write
$
(i,S)=1
$
when $(i, p_j)=1$ for all$ j\le \kappa$.
We have the following partition of $\N^*$:
\begin{equation} \label{eq-disj}
\N^* = \bigsqcup_{(i,S)=1} iS.
\end{equation}
For each element $x =
{(x_k)}_{k=1}^\infty $,   $x|_{iS}$ denotes the restriction $
x|_{iS}:={(x_{i\ell_k})}_{k=1}^\infty
$, which is also viewed as an element of $\Sigma_m$.
Given a closed subset $\Omega\subset \Sigma_m$, we define a new subset of $\Sigma_m$:
\begin{equation} \label{def-Xom}
X_\Om^{(S)}:= \Bigl\{x =
{(x_k)}_{k=1}^\infty \in \Sigma_m:\ x|_{iS} \in \Om\ \ \mbox{for all}\ i,\ (i,S)=1\Bigr\}.
\end{equation}

\begin{theorem} \label{thPSSS}
{\rm(\cite{PSSS})} There exists a vector
$
\overline{t}=(t(u))_{u\in \Pref(\Om)}\in [1,+\infty)^{\Pref(\Om)}
$
defined on the tree of prefixes of $\Omega$
such that
\begin{equation*}
t(\varnothing) \in [1,m],\ \ t(u) \in [1, m^{\ell_k(\ell_{k+1}^{-1} +\ell_{k+2}^{-1}+\cdots)}],\ |u|=k,\ k\ge 1,
\end{equation*}
which is the solution of the system
\begin{eqnarray*}
t(\varnothing)^{\gam(S)}  & =  & \sum_{j=0}^{m-1} t(j),\\
 t(u)^{\ell_{k+1}/\ell_{k}}  & = &  \sum_{j:\ uj\in \Pref_{k+1}(\Om)} t(uj),\ \ \ \forall\ u \in \Pref_k(\Om),\ \forall k\ge 1.
\end{eqnarray*}
The  Hausdorff dimension and the box dimension of $X_\Omega^{(S)}$
are respectively equal to
 \begin{eqnarray*}
\dim_H(X_\Om^{(S)}) &  =  & \log_m t(\varnothing)\\
\dim_B(X_\Om^{(S)})  & = & \gam(S)^{-1} \sum_{k=1}^\infty \Bigl(\frac{1}{\ell_k} - \frac{1}{\ell_{k+1}}\Bigr) \log_m |\Pref_{k}(\Om)|.
\end{eqnarray*}
We have $\dim_H(X_\Om^{(S)}) = \dim_B(X_\Om^{(S)})$ if and only if the tree of prefixes of $\Om$ is spherically
symmetric.
\end{theorem}

Ban, Hu and Lin \cite{BHL} studied the Minkowski dimension of $X_{2,3}$ and of some other
multiplicative sets as pattern generating problem.



\section{Remarks and Problems}

\subsection{Vector valued potential} The non-linear thermodynamic formalism can be generalized to
vectorial potentials.
 Let $\varphi, \gamma$ be two functions   defined on $\mathcal{A}^\ell$ taking real values. Instead of considering the transfer operator $\mathcal{L}_s$ as defined in (\ref{transer-operator}), we consider the following one:
$$
\mathcal{L}_s \psi (a)
= \sum_{j \in S} e^{s \varphi(a, j)+\gamma(a, j)}
\psi (Ta, j),\ a\in S^{\ell-1},\ s\in \R.
$$
There exists a unique solution to the equation
$$(\mathcal{L}_s \psi)^{\frac{1}{q}}=\psi.$$
Then we  define the pressure function $P_{\varphi,\gamma}(s)$ as indicated in
 $P_{\varphi,\gamma}(s)$.
The function $s\mapsto P_{\varphi,\gamma}(s)$ is convex and  analytic.
Now, let $\underline{\varphi}=(\varphi_1,\cdots,\varphi_d)$ be a function defined on $S^\ell$ taking values in $\R^d$. For $\underline{s}=(s_1,\cdots,s_d)\in \R^d$, we consider the following transfer operator.
$$
\mathcal{L}_{\underline{s}} \psi (a)
= \sum_{j \in S} e^{\langle \underline{s},\underline{\varphi}\rangle}
\psi (Ta, j),\ a\in S^{\ell-1},
$$
where  $\langle \cdot,\cdot\rangle$ denotes the scalar product in $\R^d$. We denote the associated pressure function by $P(\underline{\varphi})(\underline{s})$.
Then for any vectors $u,v\in \R^d$ the function
$$\R \ni s\ \longmapsto \ P(\underline{\varphi})(us+v)$$ is analytic and convex. We deduce from this that the function $\underline{s}\ \mapsto \ P(\underline{\varphi})(\underline{s})$ is infinitely differentiable and convex on $\R^d$.

Similarly, we define  the level sets $E(\underline{\alpha})$ $(\underline{\alpha}\in \R^d)$ of $\underline{\varphi}$. A vector version of Theorem \ref{thm principal} is stated by just replacing the  derivative of the pressure function by the gradient.

\subsection{Invariant spectrum and mixing spectrum}\label{InvSpectrum}
The set $E_\Phi(\alpha)$ defined by (\ref{def_ergaverage}) is not invariant.
The size of the invariant part of  $E_\Phi(\alpha)$ could be
considered to be
$$
    d_{\rm inv}(\alpha) = \sup\{\dim^* \mu:  \mu \ {\rm invariant}, \ \mu(E_\Phi(\alpha))=1\}.
$$
The function $\alpha \mapsto d_{\rm inv}(\alpha)$ is called the invariant spectrum of $\Phi$. Similarly
we define the mixing spectrum of $\Phi$ by
$$
    d_{\rm mix}(\alpha) = \sup\{\dim^* \mu:  \mu \ {\rm mixing}, \ \mu(E_\Phi(\alpha))=1\}.
$$
Examples in \cite{FSW} show that it is possible to have
$$
    d_{\rm mix}(\alpha) < d_{\rm inv}(\alpha) <d_H(\alpha).
$$

\subsection{Semigroups}
The semigroup $\{q^n\}_{n\ge 0}$ of $\mathbb{N}^*$ appeared in \cite{FSW,KPS}.
Other semigroup structures appeared in \cite{PSSS}.
Combining the ideas in \cite{PSSS,FSW}, averages like
$$
     \lim_{n \to \infty} \frac{1}{n} \sum_{k=1}^n \varphi(x_k, x_{2k}, x_{3k})
$$
can be treated \cite{Wu}. The Riesz product method
used in \cite{FLM} is well adapted to the study of the special limit on $\Sigma_2$:
$$
     \lim_{n \to \infty} \frac{1}{n} \sum_{k=1}^n (2x_k-1) (2 x_{2k} -1)\cdots (2 x_{\ell k}-1)
$$
where $\ell \ge 2$ is any integer.

\subsection{Subshifts of finite type}
What we have presented is strictly restricted to the full shift dynamics. It is a challenging problem
to study the dynamics of subshift of finite type and the dynamics with Markov property.
New ideas are needed to deal with these dynamics. It is also a challenging problem to deal with
potential depending more than one coordinates.

The doubling dynamics $Tx = 2 x$ $\!\!\!\!\mod 1$ on the interval $[0,1)$ is essentially a shift dynamics.
Cookie cutters are the first interval maps coming into the mind after the doubling map. If the cookie cutter maps are not linear, it is a difficult problem. A cookie-cutter can be coded, but the non-linearity means that the derivative is a potential depending more than one codes.

Based on the computation made in \cite{PS2}, Liao and Rams \cite{LR} considered a special piecewise linear map of two branches defined on two intervals
$I_0$ and $I_1$ and studied the following  limit
$$
          \lim_{n \to \infty} \frac{1}{n} \sum_{k=1}^n 1_{I_1}(T^kx)1_{I_1}(T^{2k}x).
$$
The techniques presented in \cite{FSW}  can be used to treat the problem
for general piecewise {\em linear} cookie cutter dynamics  \cite{FLW, Wu}.

\subsection{Discontinuity of spectrum for V-statistics}
The limit of V-statistics
$$
\lim_{n\to \infty} n^{-r}\sum_{1\le i_1, \cdots, i_r\le n} \Phi(T^{i_1}x, \cdots, T^{i_r} x).
$$
was studied in \cite{FSW_V} where it is proved that
the multifractal spectrum of topological entropy of the above limit is expressed by an variational
principle when the system satisfies the specification property.  Unlike the classical case ($r=1$)
where the spectrum is an analytic function when $\Phi$ is H\"{o}lder continuous, the spectrum of the limit of  higher order
 V-statistics ($r\ge 2$) may be  discontinuous even for very regular kernel $\Phi$.
 It is an interesting problem to determine the number of discontinuities. M. Rauch \cite{Rauch}
 has recently established a variational principe relative to $V$-statistics.



\subsection{Mutual absolute continuity of two Riesz products}

Let us state two conjectures.  See \cite{BFP} 
for the discussion on these conjectures.
\medskip

{\em Conjecture 1.}  {\em Let $\mu_a$ and $\mu_b$ be two Riesz products
and let $\omega :=(\omega_n) \in G^\mathbb{N}$. Then $\mu_a \ll
  \mu_b \Rightarrow \mu_{a,\omega} \ll \mu_{b,\omega}$, and $\mu_a
  \perp \mu_b \Rightarrow \mu_{a,\omega} \perp \mu_{b,\omega}$}.
\medskip

For a function $f$ defined on $G$, we use $\mathbb{E} f$ to denote
the integral of $f$ with respect to the Haar measure. The truthfulness is that
the preceding conjecture implies  the following one. \medskip

{\em Conjecture 2.}  {\em Let $\mu_a$ and $\mu_b$ be two Riesz products. Then\\
$$
\prod_{n=1}^\infty \mathbb{E} \sqrt{(1 + \mbox{\rm Re}\  a_n \gamma_n
 ) (1 + \mbox{\rm Re} \ b_n \gamma_n ) } >0 \Longrightarrow \mu_a \ll \mu_b;$$
$$
\prod_{n=1}^\infty \mathbb{E} \sqrt{(1 + \mbox{\rm Re}\  a_n \gamma_n
 ) (1 + \mbox{\rm Re} \ b_n \gamma_n ) } =0 \Longrightarrow \mu_a \perp
\mu_b.$$}

\subsection{Doubling and tripling}

For any integer $m\ge 2$, we define the dynamics $\tau_m x = mx $ ($\!\!\!\!\mod 1$)
on $[0, 1)$.
A typical couple of commuting transformations is the
couple $(\tau_2, \tau_3)$.
Let us take, for example, $\Phi(x, y) = e^{2\pi i (a x + b y)}$
with $a, b$ being two fixed integers. We are then led to the multiple ergodic averages,
a special case of \eqref{def_generalergaverage},
\begin{equation}\label{2_3}
   A_n^{(2, 3)}(x): = \frac{1}{n}\sum_{k=1}^n e^{2\pi i (a 2^k + b 3^k) x}.
\end{equation}
This is an object not yet well studied in the literature (but if $a=0$, we get a classical Birkhoff average). We propose to develop a thermodynamic formalism by studying Gibbs type
measures which are weak limits $\mu_{s, t}$ ($s, t \in\R$) of
$$
    Z_n(s, t)^{-1} Q_n(x) d x
$$
where
$$
     Q_n(x) := \prod_{k=1}^n e^{s \cos (2\pi (a 2^k + b 3^k) x) + t \sin (2\pi (a 2^k + b 3^k) x)}.
$$
The pressure
 function defined by
 $$
    P(s, t) := \lim_{n\to \infty} \frac{\log Z_n(s, t)}{n}
 $$
 would be differentiable. But first we have to prove the existence of the limit defining $P(s,t)$. 

 More generally, let $(c_n)$ be a sequence of complex numbers and $(\lambda_n)$
 a lacunary sequence of positive integers (by lacunary we mean $\inf_n \frac{\lambda_{n+1}}{\lambda_n}>1$). We can consider the following weighted
 lacunary trigonometric averages
 $$
    \frac{1}{n}\sum_{k=1}^n c_k e^{2\pi i \lambda_k x}.
 $$
 Under the divisibility condition $\lambda_n | \lambda_{n+1}$, such averages and more
 general averages were studied in \cite{Fan1997b}. For example, if $c_k = e^{2\pi i \omega_k}$
 with $(\omega_k)$ being an i.i.d. sequence of Lebesgue distributed random variables,  from the results obtained in \cite{Fan1997b} we deduce that almost surely the pressure
 is well defined and equal to the following deterministic function
 $$
     P(s, t) = \log \int_0^{2\pi} e^{\sqrt{t^2+s^2} \cos x}\frac{ d x}{2\pi }.
 $$
 Recall that $$
 J_0(r) = \frac{1}{2\pi}\int_{0}^{2\pi} e^{r\cos x} d x
 = \sum_{n=0}^\infty \frac{r^{2 n}}{(n!)^2 2^{2n}}.
  $$
  is the Bessel function.

 But,  the condition $\lambda_n |\lambda_{n+1}$ for $\lambda_n = 2^n + 3^n$
  is not satisfied. Neither the condition is  satisfied for $\lambda_n =2^n + 4^n$.
  No rigorous results are known for the multifractal analysis of the averages defined by (\ref{2_3}).

 As conjectured by
F\"urstenberg, the Lebesgue measure is the unique continuous probability measure which is both $\tau_2$-invariant and $\tau_3$-invariant.
However,
common $\tau_2$- and $\tau_3$-periodic points (different from the trivial one $0$) do exists. Given two integers $n\ge 1$ and $m\ge 1$.
We can prove that there is a point $x (\not=0)$ which is $n$-periodic with respect to $\tau_2$
and $m$-periodic with respect to $\tau_3$ if and only if
\begin{equation}\label{2_3condition}
    (2^n -1, 3^m-1) >1.
\end{equation}
Let $d = (2^n -1, 3^m-1)$. When the above condition on GCD is satisfied, there are $d-1$ such common periodic
points. These common periodic points $x (\not=0)$ are of the form
$
    x = \frac{k}{2^n-1}= \frac{j}{3^m-1}
    $
    for some
    $1\le k < 2^n-2$ and  $1\le j <3^m-2$.
    Actually choices for $k$ are
$$
            1 \cdot \frac{2^n-1}{d},\ \ \  2 \cdot \frac{2^n-1}{d},\ \cdots, \  (d-1) \cdot \frac{2^n-1}{d}.
           $$
           Choices for $j$ are
$
            1 \cdot \frac{3^m-1}{d},\ \ \  2 \cdot \frac{3^m-1}{d},\ \cdots, \  (d-1) \cdot \frac{3^m-1}{d}.
           $
Thus the $d-1$ common periodic points are $\frac{1}{d}, \frac{2}{d}, \cdots, \frac{d-1}{d}$.
For such a point $x$, the following limit exists
$$
   \lim_{N\to \infty} A_N^{(2,3)}(x) = A_{nm}^{(2,3)}(x).
$$
Note that there is an infinite number of such couples $n$ and $m$ such that (\ref{2_3condition}) holds.
There would be some relation between these common periodic points and the multifractal behavior of
$A_N^{(2,3)}(x)$.



\def\refname{Bibliography}

\end{document}